\documentclass[12pt]{amsart}%
\usepackage{graphicx}
\usepackage{amssymb}
\usepackage{amsmath}
\usepackage{amsfonts}
\numberwithin{equation}{section}
\begin{document}
\title[Thomas-Fermi Equation]{Dynamically Consistent Approximate Rational Solutions to the Thomas-Fermi
Equation }
\author[R. E. Mickens and I. Herron]{Ronald E. Mickens\\
Department of Physics \\
Clark Atlanta University \\
Atlanta, GA 30314 USA \\
Email: rmickens@cau.edu \\
 and \\
Isom H. Herron \\
Department of Mathematical Sciences\\
Rensselaer Polytechnic Institute\\
Troy, NY 12180 USA \\
Email: herroi@rpi.edu}

\begin{abstract}
We construct two rational approximate solutions to the Thomas-Fermi (TF)
nonlinear differential equation. These expressions follow from an application
of the principle of dynamic consistency. In addition to examining differences
in the predicted numerical values of the two approximate solutions, we compare
these values with an accurate numerical solution obtained using a fourth-order
Runge-Kutta method. We also present several new integral relations satisfied
by the bounded solutions of the TF equation

\end{abstract}

\maketitle

\section{Introduction}

\label{sec;introduction}

The major purpose of this paper is to construct a new rational approximation
to the Thomas-Fermi equation
\begin{equation}
\frac{d^{2}y\left(  x\right)  }{dx^{2}}=\frac{y\left(  x\right)  ^{\frac{3}%
{2}}}{x^{\frac{1}{2}}}, \label{1.1}%
\end{equation}
with boundary conditions
\begin{equation}
y\left(  0\right)  =1,y\left(  \infty\right)  =0. \label{1.2}%
\end{equation}

This equation arises in the modeling of certain phenomena in the atomic
structure of matter (\cite{1}, \cite{2}).

Since it is expected that no exact analytical solution formula exists for the
general solution in terms of a finite combination of elementary functions, a
broad range of approximate solutions have been derived (\cite{3}, \cite{4},
\cite{5}, \cite{6}, \cite{7}, \cite{8}). The mathematical procedures used to
obtain these expressions include variational methods (\cite{3}, \cite{6}), the
use of iteration techniques \cite{4}, application of the Adomian decomposition
method \cite{5}, homotopy analysis \cite{7}, and rational approximations
\cite{8}. The indicated references (\cite{3}, \cite{4}, \cite{5}, \cite{6},
\cite{7}, \cite{8}) are but a small subset of the hundreds of publications on
this subject. Note that the existence of such a substantial research
literature is due, in large part, to the fact that this thorny boundary-value
problem has drawn out a significant number of ideas as to how approximations
to its solutions should be done.

The purpose of this paper is to use a new methodology to construct approximate
solutions to the Thomas-Fermi (TF) equation. We do this within the framework
of the concept or principle of \textit{dynamic consistency} \cite{9} and the
restrictions that it imposes on any mathematical representation of the forms
used for rational approximations to a TF equation solution. We demonstrate
that this is rather easy to do.

In outline, Section 2 provides a summary of several of the general, exact, and
for the most part qualitative properties of the TF equation \cite{10}. This
section includes several new results of the authors. Section 3 introduces and
defines the concept of \textit{dynamic consistency} \cite{9} and it is used in
Section 4 to construct two elementary rational approximations to the TF
solution. Finally, in the last section, we provide a summary of our results
and discuss possible extensions to this work.

\section{Preliminaries: Exact results}

\label{sec;preliminary}

It turns out that even in the absence of knowledge of an explicit, exact
solution for the TF equation many of the essential properties of such a
solution can be derive. Below, we give a concise summary of these items and
refer the reader to the paper of Hille \cite{10} for proofs of some of these statements.

a) An exact solution to the TF equation is
\begin{equation}
y\left(  x\right)  \equiv y_{s}\left(  x\right)  =\frac{144}{x^{3}}.
\label{2.1}%
\end{equation}
Note that it contains no arbitrary parameters and thus is not a special case
of the (unknown) general solution. We denote it as, $y_{s}\left(  x\right)  $
and call it a singular or asymptotic solution in the sense that given a
solution $y\left(  x\right)  $ then \cite{10},
\begin{equation}
\lim_{x\rightarrow\infty}(x^{3}y\left(  x\right)  )=144. \label{2.2}%
\end{equation}

b) The curve $y\left(  x\right)  \equiv y_{s}\left(  x\right)  \ $separates
the bounded and unbounded solutions of the TF equation. This is evident from
inspection of the flow-space as indicated in Figure 1.%

\begin{figure}
[ptb]
\begin{center}
\includegraphics[
height=1.9993in,
width=3.0007in
]%
{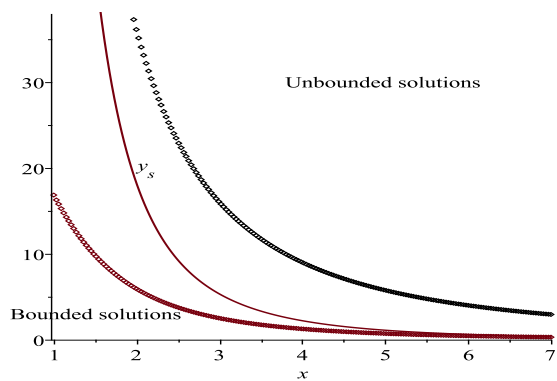}%
\caption{Singular solution $y_{s}$ separating unbounded and bounded solutions}%
\end{center}
\end{figure}

Note that since $y\left(  x\right)  \equiv0$, is a solution, it follows that
all solutions below the $y=y_{s}\left(  x\right)  $ curve are bounded and all
solutions above are unbounded. The bounded solutions have a finite value for
$y\left(  0\right)  , $ while the unbounded have $y\left(  0\right)  =+\infty$.

c) The physical phenomena for which the TF equation was constructed requires
\begin{align}
y\left(  x\right)   &  >0,y^{\prime}\left(  x\right)  \equiv\frac{dy\left(
x\right)  }{dx}<0,0\leq x\,<\infty,\label{2.3a}\\
y\left(  0\right)   &  =1,\ y\left(  \infty\right)  =0. \label{2.3b}%
\end{align}
Observe that since $y^{\prime\prime}\left(  x\right)  >0,0<x<\infty,$ then all
solutions are concave upward and this is consistent with the trajectory-flow
given in Figure 1.

The figure also provides an overview of the bounded and unbounded trajectories
of the TF equation.

d) For the bounded solutions%
\begin{equation}
0<y\left(  x\right)  <y\left(  0\right)  <\infty, y^{\prime}\left(  0\right)
<y^{\prime}\left(  x\right)  <0,y^{\prime\prime}\left(  x\right)
>0,0<x\,<\infty. \label{2.4}%
\end{equation}
Further we have
\begin{equation}
y^{\prime\prime}\left(  x\right)  =O\left(  \frac{1}{\sqrt{x}}\right)
,x\rightarrow0^{+}. \label{2.5}%
\end{equation}

e) It is remarkable that $y^{\prime}\left(  0\right)  $ can be calculated to
essentially any decimal place of accuracy \cite{11}; defining the constant $B$
as follows,
\begin{equation}
y^{\prime}\left(  0\right)  =-B, \label{2.6}%
\end{equation}
we have \cite{11}
\begin{equation}
B=1.588071022611375\ldots. \label{2.7}%
\end{equation}

f) For small $x$, $y\left(  x\right)  $ has the representation \cite{4} for
$y\left(  0\right)  =1$ and $y^{\prime}\left(  0\right)  =-B,$
\begin{align}
y\left(  x\right)   &  =\left[  1-Bx+\left(  \frac{1}{3}\right)  x^{3}-\left(
\frac{2B}{15}\right)  x^{4}+\cdots\right]  \nonumber\\
&  \cdots+x^{3/2}\left[  \frac{4}{3}-\left(  \frac{2B}{5}\right)  x+\left(
\frac{3B^{2}}{70}\right)  x^{2}\right.  +\left.  \left(  \frac{2}{27}%
+\frac{B^{3}}{252}\right)  x^{3}+\cdots\right]  .\label{2.8}%
\end{align}

g) The following \textquotedblleft sum-rules\textquotedblright, involving
integrals of $y\left(  x\right)  ,$have been derived by the authors. As far as
we are aware, most of these are \ new results pertaining to the solution (see
Eqs. (\ref{1.1}) and (\ref{1.2})) of the TF equation:%
\begin{equation}
\int_{0}^{\infty}\sqrt{x}y\left(  x\right)  ^{3/2}dx=1, \label{2.9}%
\end{equation}
(was derived in \cite{15})%
\begin{equation}
\int_{0}^{\infty}y\left(  x\right)  dx=\frac{1}{2}\int_{0}^{\infty}x^{\frac
{3}{2}}y\left(  x\right)  ^{\frac{3}{2}}dx, \label{2.10}%
\end{equation}%
\begin{equation}
B=\int_{0}^{\infty}\left[  \left(  \frac{dy}{dx}\right)  ^{2}+\frac
{y^{\frac{5}{2}}}{\sqrt{x}}\right]  dx \label{2.11}%
\end{equation}%
\begin{equation}
B=\int_{0}^{\infty}\frac{y\left(  x\right)  ^{3/2}}{\sqrt{x}}dx. \label{2.12}%
\end{equation}

\bigskip Let us derive the result in Eq. (\ref{2.11}). Multiply both sides of
the TF equation by $y\left(  x\right)  :$
\begin{equation}
yy^{\prime\prime}=\frac{y^{\frac{5}{2}}}{\sqrt{x}}. \label{2.13}%
\end{equation}
Integrating once by parts and evaluating at $x=0$ and $x=\infty,$ gives%
\begin{equation}
y\left(  \infty\right)  y^{\prime}\left(  \infty\right)  -y\left(  0\right)
y^{\prime}\left(  0\right)  -\int_{0}^{\infty}\left(  y^{\prime}\left(
x\right)  \right)  ^{2}dx=\int_{0}^{\infty}\frac{y^{\frac{5}{2}}}{\sqrt{x}}dx.
\label{2.14}%
\end{equation}
But $y\left(  0\right)  =1,y^{\prime}\left(  0\right)  =-B$ and $y\left(
\infty\right)  =y^{\prime}\left(  \infty\right)  =0;$ therefore, the result
listed in Eq. (\ref{2.11}) is obtained. The other relations can be derived
using similar techniques.

Finally, we go back to Eq. (\ref{2.8}) and demonstrate one possible way to
derive it. First, observe that for small $x,$ the TF equation is approximated
by the expression%
\begin{equation}
y^{\prime\prime}\left(  x\right)  \simeq\frac{1}{\sqrt{x}}, \label{2.15}%
\end{equation}
where $y\left(  0\right)  =1$ and $y^{\prime}\left(  0\right)  =-B.$
Integrating twice and imposing the just stated initial conditions gives%
\begin{equation}
y\left(  x\right)  \simeq1-Bx+\left(  \frac{4}{3}\right)  x^{\frac{3}{2}},
\label{2.16}%
\end{equation}
corresponding to the first three terms of (\ref{2.8}).

Let us now construct an iteration scheme to obtain a small$\ x$ representation
for the TF equation. We assume it has the form%
\begin{align}
y_{N+1}^{\prime\prime}(x)  &  =\frac{\left[  y_{N}\left(  x\right)  \right]
^{\frac{3}{2}}}{\sqrt{x}},\ N=0,1,2,\ldots;\label{2.17a}\\
y_{N}\left(  0\right)   &  =1,y_{N}^{\prime}\left(  0\right)
=-B,\ \label{2,17b}\\
y_{0}\left(  x\right)   &  =1. \label{2.17c}%
\end{align}
Note that $y_{1}\left(  x\right)  $ is the result given in Eq. (\ref{2.16}).

To calculate $y_{2}\left(  x\right)  $, we start with%
\begin{align}
y_{2}^{\prime\prime}(x)  &  =\frac{\left[  y_{1}\left(  x\right)  \right]
^{\frac{3}{2}}}{\sqrt{x}}\nonumber\\
&  =\frac{\left[  1-Bx+\left(  \frac{4}{3}\right)  x^{\frac{3}{2}}\right]
^{\frac{3}{2}}}{\sqrt{x}}. \label{2.18}%
\end{align}
Expanding the expression in the brackets and retaining the first two terms
gives%
\begin{equation}
\left[  1-Bx+\left(  \frac{4}{3}\right)  x^{\frac{3}{2}}\right]  ^{\frac{3}%
{2}}=1-\left(  \frac{3B}{2}\right)  x+\cdots\ . \label{2.19}%
\end{equation}
Thus it follows that
\begin{equation}
y_{2}^{\prime\prime}\left(  x\right)  =\frac{1}{\sqrt{x}}-\left(  \frac{3B}%
{2}\right)  \sqrt{x}. \label{2.20}%
\end{equation}
Integrating this equation twice, with the imposition of the initial conditions
$y_{2}\left(  0\right)  =1,y_{2}^{^{\prime}}\left(  0\right)  =-B,$ gives%
\begin{equation}
y\left(  x\right)  =1-Bx+\left(  \frac{4}{3}\right)  x^{\frac{3}{2}}-\left(
\frac{2B}{5}\right)  x^{\frac{5}{2}}. \label{2.21}%
\end{equation}
The result corresponds, respectively, to the first two terms in the brackets
of Eq. (\ref{2.18}).

Finally, with careful attention to which terms to retain, at each step of the
iteration given by Eq. (\ref{2.17a}-\ref{2.17c}), the full result presented in
Eq. (\ref{2.8}) can be produced.

\section{Dynamic consistency}

Consider two systems $A$ and $B.$ Let system $A$ have the property $P.$ If
system $B$ also has property $P,$ then $B$ is said to be dynamic consistent to
$A$, with respect to property $P\ $\cite{9}.

Dynamic consistency \ (DC) has served as one of the fundamental principles of
the nonstandard, finite difference methodology for constructing improved
discrete models for the numerical integration of differential equations
\cite{9}. It serves as an assessment of the \textquotedblleft
closeness\textquotedblright\ or \textquotedblleft fidelity\textquotedblright%
\ of two systems based on their mutually shared properties. Note that the
systems $A$ and $B$ may not be of the same type object, structure or thing.
For example, $A$ might be a physical system, while $B$ could be a mathematical
model of $A.\ $Likewise, $A$ might be a system of differential equations and
$B$ a particular discretization of these equations. Thus DC is of general
applicability and is available for use in diverse situations.

In the next section, we provide two ansatzes for approximate solutions to the
TF equation, Our proposed solutions are to be DC with respect to the following
properties of the exact, but unknown, solution to the second-order, nonlinear
TF differential equation:%
\begin{equation}
\text{(i)}\ \ y\left(  0\right)  =1,y^{\prime}\left(  0\right)  =-B;
\label{3.1}%
\end{equation}%
\begin{equation}
\text{(ii)}\ \ 0<y\left(  0\right)  \leq1,-B\leq y^{\prime}\left(  x\right)
<0,0\leq x<\infty; \label{3.2}%
\end{equation}%
\begin{equation}
\text{(iii)}\ \ \lim_{x\rightarrow\infty}\left(  x^{3}y\left(  x\right)
\right)  =144\ ;\ \label{3.3}%
\end{equation}

(iv) the approximate solutions $y_{a}\left(  x\right)  $ are to be taken as
rational functions of either $x$ or $\sqrt{x}$;

(v) the approximate solutions$, y_{a}\left(  x\right) $, are to have one or
the other of the following forms for small $x$;%

\begin{equation}
y_{a}\left(  x\right)  =\left\{
\begin{array}
[c]{c}%
1-Bx,\ \ \text{rational in }x;\\
1-Bx+\left(  \frac{4}{3}\right)  x^{3/2},\text{ rational in } \sqrt{x}%
\end{array}
\right.  \label{3.4}%
\end{equation}

\section{Approximate solutions}

\label{sec;distribution}

We take the following two ansatzes to represent approximations to the solution
of the TF equation
\begin{equation}
y_{a}^{\left(  1\right)  }(x)=\frac{1}{1+Bx+\left(  \frac{1}{144}\right)
x^{3}}, \label{4.1}%
\end{equation}%
\begin{align}
y_{a}^{\left(  2\right)  }\left(  x\right)   &  =\frac{1}{1+Bx-\left(
\frac{4}{3}\right)  x^{3/2}+Cx^{2}+\left(  \frac{1}{144}\right)  x^{3}%
},\label{4.2a}\\
C  &  =\frac{B^{2}}{2}, \label{4.2b}%
\end{align}
where $B=1.588071$, and as a consequence of (\ref{4.2b}), $C=1.260985.$

Close inspection of $\ y_{a}^{\left(  1\right)  }(x)$ and $y_{a}^{\left(
2\right)  }\left(  x\right)  $ shows that they are DC with the five conditions
listed in section 3. Observe that $y_{a}^{\left(  1\right)  }\left(  x\right)
$ is a rational function of $x$, while $y_{a}^{\left(  2\right)  }\left(
x\right)  $ is a rational function of $\sqrt{x}.$%

\begin{figure}
[ptb]
\begin{center}
\includegraphics[
height=2in,
width=3in
]%
{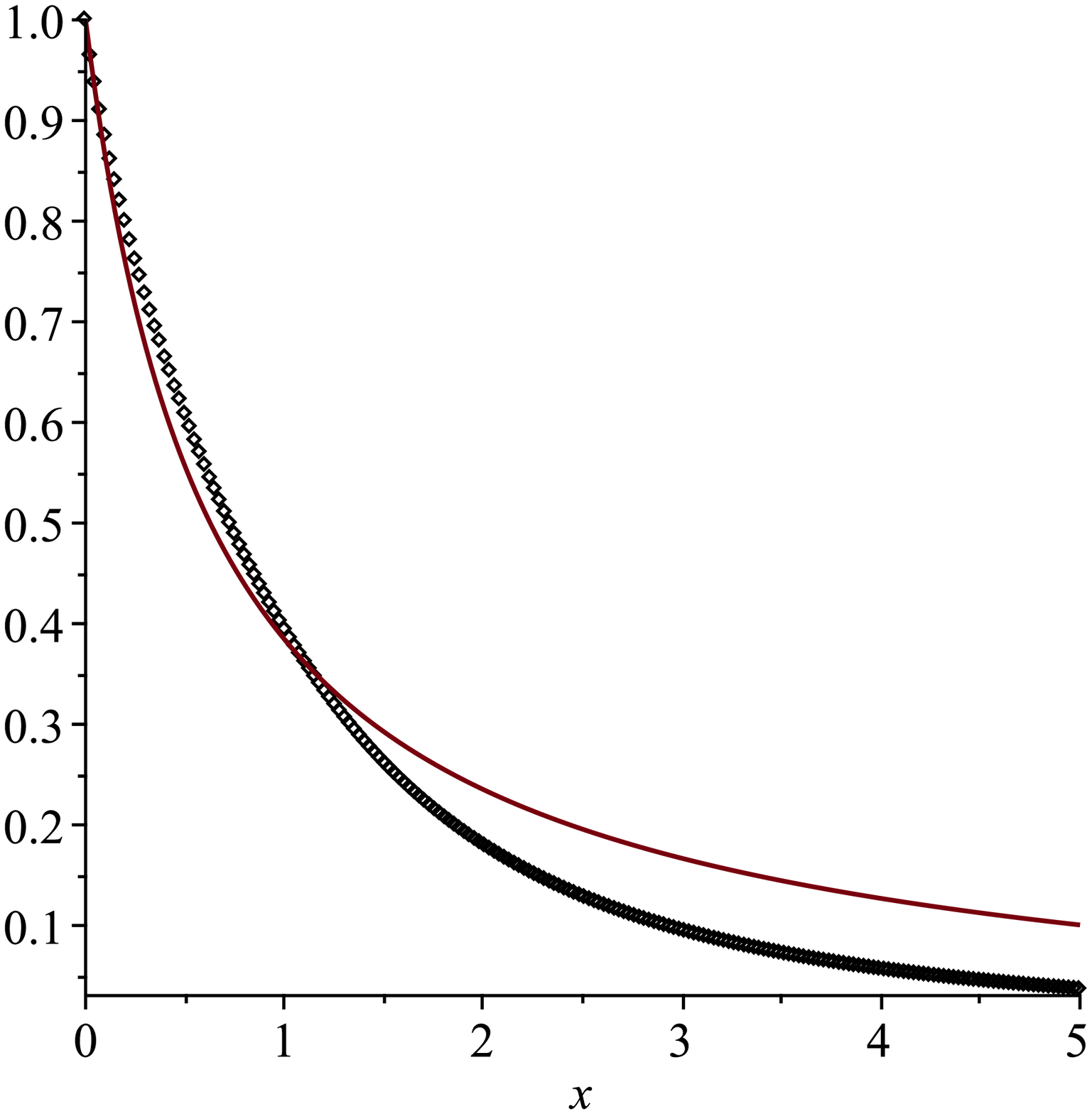}%
\caption{Plots of $y_{a}^{\left(  1\right)  }\left(  x\right)  $ $\left(
^{\_\_\_\_\_}\right)  \ $and $y_{a}^{\left(  2\right)  }\left(  x\right)
\ \left(  \cdot\cdot\cdot\cdot\right)  $} 
\end{center}
\end{figure}

Figure \ref{fig:2} gives plots of $y_{a}^{\left(  1\right)  }(x)$ and
$y_{a}^{\left(  2\right)  }(x)$ vs. $x$. The two curves lie close to each
other, but do not coincide, i.e.%
\begin{equation}
\left\{
\begin{array}
[c]{c}%
y_{a}^{\left(  1\right)  }(x)<y_{a}^{\left(  2\right)  }(x),\ 0<x<x_{0},\\
y_{a}^{\left(  1\right)  }(x)>y_{a}^{\left(  2\right)  }(x),\ \ x>x_{0}%
\end{array}
\right.  \label{4.3}%
\end{equation}
with $x_{0}\approx1.$

\section{Discussion}

\label{sec5}

In summary we have constructed two representations of approximate solutions to
the TF equation. However, it must be stated that this equation is itself an
approximation to an (unknown) equation describing atomic quantum phenomena
(\cite{1}, \cite{2}). In addition, all numerical procedures have errors
associated with their algorithms and the hardware on which they are implemented.

The methodology used to generate the ansatzes, given in Eqs. (\ref{4.1}) and
(\ref{4.2a}, \ref{4.2b}), is based on constructing for $y\left(  x\right)  $,
the exact solution, rational functions such that the five conditions listed in
section 3 are satisfied. An important benefit of this procedure is that the
process is essentially algebraic and requires the introduction of no
parameters which must be determined later in the construct. The $y_{a}%
^{\left(  1\right)  }(x)$ and $y_{a}^{\left(  2\right)  }(x)$ are constructed
to be consistent with the $x-$small and $x-$large asymptotics. Thus, these
representations possess great simplicity both in their manner of construction,
as well as the lack of a need to introduce new parameters.

In Table 1 we compare the numerics of a fourth-order Runge-Kutta method and
the values of $y\left(  x\right)  $ as calculated from our formulae for
$y_{a}^{\left(  1\right)  }(x)$ and $y_{a}^{\left(  2\right)  }(x)$. Clearly,
both analytic approximations and the numerical integration results are
consistent, and in fact are in very good qualitative agreement with each other.

\begin{center}
\bigskip Table - 1 Numerical values of $y\left(  x\right)  $ and $y_{a}\left(
x\right)  .$ \newline%
\begin{tabular}
[c]{lll}%
$x$ & $y\left(  x\right)  $ numerical \cite{12} & $y_{a}^{\left(  1\right)
}\left(  x\right)  $\\
$0.0$ & $1.0000$ & $1.0000$\\
$0.5$ & $0.6070$ & $0.5571$\\
$1.0$ & $0.4240$ & $0.3853$\\
$2.0$ & $0.2430$ & $0.2363$\\
$4.0$ & $0.1840$ & $0.1283$\\
$5.0$ & $0.0789$ & $0.1020$\\
$10.0$ & $0.0243$ & $0.0420$\\
$25.0$ & $0.0035$ & $0.0067$\\
$40.0$ & $0.0011$ & $0.0020$%
\end{tabular}%
\begin{tabular}
[c]{l}%
$y_{a}^{\left(  2\right)  }\left(  x\right)  $\\
$1.0000$\\
$0.6102$\\
$0.3964$\\
$0.1817$\\
$0.0578$\\
$0.0378$\\
$0.0093$\\
$0.0013$\\
$0.0005$%
\end{tabular}

\end{center}

Since the main task has been achieved, i.e., to demonstrate through explicit
construction an alternative methodology for finding rational approximations
for the solution of the TF equation, we now can proceed to the next stage,
namely, determine a restriction or condition allowing us to access the overall
validity or accuracy of our solutions. To do this, we consider the integral
relations stated in Section 2. For our purposes, we select the sum-rule given
in Eq. (\ref{2.11}), mainly because it takes the form of an energy type
integral \cite{13}, widely used throughout physics. In fact, inspection of Eq.
(\ref{2.11}) suggests the following correspondence%
\[
\left(  \frac{dy}{dx}\right)  ^{2}\longrightarrow\text{ kinetic energy,}%
\]%
\[
\frac{y^{\frac{5}{2}}}{\sqrt{x}}\longrightarrow\ \text{potential energy.}%
\]
Note that the result in Eq. (\ref{2.11}) can be used as a consistency check,
since $B$ occurs in $y\left(  x\right)  $ which appears in the integral. Using
$y_{a}^{\left(  1\right)  }(x)$ and $y_{a}^{\left(  2\right)  }(x)$ and
carrying out the required numerical integration we find, respectively%
\[
B_{1}=1.584744,\ \ B_{2}=1.592931,
\]
which are to be compared to the input value of%
\[
B=1.588071.
\]
The fractional percentage errors are%
\[
E_{1}=\left(  \frac{B-B_{1}}{B}\right)  \cdot100=0.21\%,\ E_{2}=\left(
\frac{B-B_{2}}{B}\right)  \cdot100=-0.27\%.
\]
Thus, with regard to the sum-rule, stated in Eq. (\ref{2.11}), both rational
approximations are essentially the same in their ability to satisfy the sum-rule.

Eventually, the next step in this research is to construct more complex
rational approximations for the solution of the TF equation which incorporate
further terms in the known asymptotic expansions for both $x-$small and
$x-$large using Pad\'{e} approximants \cite{14}.




\section*{}


\begin{thebibliography}{99}                                                                                               %


\bibitem {1}E. Fermi,
\textit{Rend. Acad. Naz. Lincei} \textbf{6} (1927), 602-607

\bibitem {2}L. H. Thomas,
\textit{Proc. Camb. Phil. Soc. }\textbf{23} (1927), 452-548.

\bibitem {3}N. Anderson, A. M. Arthurs, and P. D. Roninson,
\textit{Il Nuovo Cimento B Series} \textit{10}, \textbf{57} (1968), 523-526.

\bibitem {4}G. I. Plindev and S. K. Pogrebnga,
\textit{Journal of Physics B: Atomic and Molecular Physics 20}, Number 17
(1987), L547.

\bibitem {5}G. Adomian,
\textit{Applied Mathematics Letters} \textbf{11}, No. 3 (1998), 131-133.

\bibitem {6}M. Desaix, D. Anderson, and M. Lisak,
\textit{European Journal of Physics} \textbf{25}, No. 6 (2004), 699-705.

\bibitem {7}(A. El-Nahhas,
\textit{Acta Physica Polonica A} \textbf{114}, No. 4 (2008), 913-918.

\bibitem {8}W. Robin,
\textit{Journal of Innovative Technology and Education} 5, No. 1 (2018), 7-13.

\bibitem {9}R. E. Mickens,
\textit{Journal of Difference equations and Applications} \textbf{11}, No.7
(2005), 645-653.

\bibitem {10}E. Hille,
\textit{Proc. Nat. Acad Sci. USA} \textbf{62}, no. 1 (1969), 7-10.

\bibitem {11}F. M Fern\'{a}ndez,
\textit{Physics Letters A} 365 (2007), 111

\bibitem {12}P. S. Lee and T.-Y. Wu,
\textit{Chinese Journal of Physics} \textbf{35}, No. 6-11 (1997), 737-741.

\bibitem {13}J. R. Taylor,\textit{ Classical Mechanics}, University Science
Books, Sausalito, CA, 2005.

\bibitem {14}M. A. Noor and S. T. Mohyud-Din,
\textit{International Journal of Nonlinear Science} \textbf{8}, No. 1 (2009), 27-31

\bibitem {15}E. B. Baker,
\textit{Physical Review} \textbf{16 }(1930), 630-647.
\end{thebibliography}
\end{document}